# Maximum distance between consecutive primes and other related questions

Victor Volfson

ABSTRACT  The paper considers the asymptotic of the ratio of the number of primes not exceeding the primorial and the number of residues in the reduced system of residues for the given primorial. We study the relationship between asymptotic lower bounds for the values of the Jacobstal function and the maximum distance between successive primes. One algorithm for computing the Jacobstal function is given. The substantiation of the conjectures about the upper estimate of the maximum distance between successive prime numbers is given.

Keywords: reduced system of residues, primorial, Jacobstal function, algorithm for calculating the Jacobstal function, connection of asymptotic estimates of the distance between residues and successive primes, Legendre's conjecture, conjecture on the upper estimate of the maximum distance between successive primes.

## 1. INTRODUCTION

The reduced residue system modulo M (RRS M) are natural numbers from 1 to M that do not have common divisors (coprime) with M. The primorial of a prime number is the product $p\# = \prod_{2 \leq q \leq p} q$, where $q$ are consecutive primes. We will consider the reduced residue system modulo equal to the primorial - RRS $p\#$.

It is known that there are removed (from the natural series) multiples of $p_r$ natural numbers, where $r$ is the ordinal number of a prime, at the $r$-th step of the sieve of Eratosthenes. Based on the properties of the sieve of Eratosthenes, all remaining natural numbers on the interval from $p_{r+1}$ to $p_{r+1}^2$ are prime after the $r$-th step of the sieve of Eratosthenes.

If we take $p = p_r$ in RRS $p\#$, then it will not contain all natural numbers that are multiples of consecutive prime numbers from 2 to $p_r$, similarly to the sieve of Eratosthenes, therefore all residues will also be prime numbers in the interval from $p_{r+1}$ to $p_{r+1}^2$. This is the property of RRS $p\#$, coinciding with the property of the sieve of Eratosthenes. We will use it in the future.

Thus, RRS $p\#$ starts with residues of 1, $p_{r+1}$ with the distance between residues much greater than between prime numbers: $p_{r+1} - 1$, then comes the interval from $p_{r+1}$ to $p_{r+1}^2$, where there are only prime numbers, then there are residues that are not prime numbers. Therefore, the question arises - what is the asymptotic of the ratio of the number of primes that do not exceed $p_r\#$, and the number of residues in RRS $p\#$ ? The answer to this question is given in chapter 2 of this work.

This issue is closely related to the issue - the distance between successive residues in RRS $p\#$ and in particular with the maximum distance between successive residues in RRS $p\#$. The arithmetic function of the dependence of the maximum distance between successive residues in RRS $p\#$ on the ordinal number of a prime number is called the Jacobstal function. The values of this function are given in table A048670 QEIS. This table shows that the Jacobstal function is strictly increasing.

The works [1], [2], [3] are devoted to algorithms for calculating the values of the Jacobstal function. Improving the efficiency of algorithms for calculating the values of the Jacobstal function is an urgent task. Chapter 3 will discuss one of these algorithms.

Asymptotic upper and lower bounds for the Jacobstal function are given in [4], [5], [6]. The last proven asymptotic lower estimate for the values of the Jacobstal function is:

$$j(x\#) \gg x \ln x \ln \ln \ln x / \ln \ln x. \qquad (1.1)$$

The last proven asymptotic upper bound for the values of the Jacobstal function:

$$a(k) \ll k^2 (\ln k)^2, \qquad (1.2)$$

where $k$ is the ordinal number of the prime.

We will talk about the connection between asymptotic lower bounds for the values of the Jacobstal function and the maximum distance between successive primes in Chapter 5 of this paper.

Let's take the value $p = p_{r+1}$ in RRS $p\#$. All residues from $p_{r+2}$ to $p_{r+2}^2$ will be prime numbers in this case. Thus, (in comparison with RRS $p_r\#$) a new interval from $p_{r+2}^2$ to $p_{r+2}^2$ will be added to RRS $p_{r+1}\#$, on which all residues will be prime numbers.

Therefore, as the value $p$ increases in RRS $p\#$ from $p = 2$ onwards, the union of such intervals will include all primes. Let's use this to estimate the maximum distance between consecutive primes.

The papers [7], [8], [9] are devoted to finding an asymptotic upper bound for the maximum distance between successive primes.

Let's denote the maximum distance between consecutive prime numbers $d = \sup_{p_{k+1} \leq x} \{p_{k+1} - p_k\}$. If we put $x = p_k$, then in these papers the asymptotic of the upper bound for the maximum distance between successive primes is represented as $d < x^\Theta$. The best proven upper asymptotic estimate at the moment is the estimate with $\Theta = 0,525$.

We extract the Legendre conjecture from unproven asymptotic upper bounds for the maximum distance between successive primes, the proof of which is one of the problems pointed out by Landau.

According to this conjecture, there is a prime number between $n^2$ and $(n+1)^2$ for any natural n.

If we put $x = n^2$, then, following the Legendre conjecture, we obtain the following requirement for the maximum distance between successive prime numbers:

$$d < (n+1)^2 - n^2 = (\sqrt{x}+1)^2 - x = 2\sqrt{x} + 1. \tag{1.3}$$

This shows that proven asymptotic estimate $d < x^{0,525}$ is inferior to estimate (1.3).

If we put $x = p_k$ in (1.3), then we get the Andrica conjecture:

$$d < 2\sqrt{p_k} + 1. \tag{1.4}$$

A stronger conjecture about an upper bound for the maximum distance between successive primes than (1.4) is substantiated in Chapter 4 of the paper.

2. ASYMPTOTIC OF THE RATIO OF THE MUMBER OF PRIMES NOT EXCEEDING $p\#$, AND THE NUMBER OF RESIDUES IN RRS $p\#$

Assertion 1

The asymptotic of the ratio of the number of primes not exceeding $p\#$, and the number of residues in RRS $p\#$ is equal to:

$$\frac{\pi(p\#)}{\varphi(p\#)} = \frac{e^\gamma \ln p}{p}(1+o(1)), \tag{2.1}$$

where $\gamma$ is the Euler constant.

Proof

Based on the properties of the Euler function:

$$\varphi(p\#) = \varphi(\prod_{q \leq p} q) = \prod_{q \leq p}(q-1) = \prod_{q \leq p} q \prod_{q \leq p}(1-1/q) = p\# \prod_{q \leq p}(1-1/q), \tag{2.2}$$

where $q$ is a prime number.

Having in mind the asymptotic law of prime numbers:

$$\pi(p\#) = \frac{p\#}{\ln(p\#)}(1+o(1)). \tag{2.3}$$

Using (2.2) and (2.3) we get:

$$\frac{\pi(p\#)}{\varphi(p\#)} = \frac{p\#}{p\# \ln(p\#) \prod_{q \leq p}(1-1/q)}(1+o(1)) = \frac{(1+o(1))}{\ln(p\#) \prod_{q \leq p}(1-1/q)}. \qquad (2.4)$$

Based on (2.4) and Mertens formula, we obtain:

$$\frac{\pi(p\#)}{\varphi(p\#)} = \frac{(1+o(1))}{\ln(p\#) \prod_{q \leq p}(1-1/q)} = \frac{e^{\gamma} \ln p}{\ln(p\#)}(1+o(1)). \qquad (2.5)$$

Considering that $\ln(p\#) = \sum_{q \leq p} \ln q = p(1+o(1))$, based on (2.5) we have:

$$\frac{\pi(p\#)}{\varphi(p\#)} = \frac{e^{\gamma} \ln p}{\ln(p\#)}(1+o(1)) = \frac{e^{\gamma} \ln p}{p}(1+o(1)),$$

which corresponds to (2.1).

### 3. ONE ALGORTHM FOR COMPUTING THE JACOBSTAL FUNCTION

Residues $p_r^2, p_r p_{r+1}, \ldots$ are removed during the transition from RRS $p_{r-1}\#$ to RRS $p_r\#$, so new distances between RRS $p_r\#$ residues appear. Since the Jacobstal function is strictly increasing, there are the largest ones among these new distances between RRS $p_r\#$ residues on the interval $[1, p_r\#]$.

Thus, it is sufficient to view only the distances between the RRS $p_r\#$ residues obtained after removing the residues $p_r^2, p_r p_{r+1}, \ldots$ to determine the value of the Jacobstal function for RRS $p_r\#$.

Having in mind the symmetry of RRS $p_r\#$, then there are an even number of such maximum distances between residues on the interval $[1, p_r\#]$. Therefore, it suffices to view only the distances between RRS $p_r\#$ residues obtained after removing residues $p_r^2, p_r p_{r+1}, \ldots$ on the interval $[1, p_r\#/2]$ to determine the Jacobstal function for RRS $p_r\#$.

For example, let's consider RRS $11\#$ on the interval $[1, 1155]$.

Removable residues: $121, 143, 187, 209, 253, 319, 341, 407, 451, 473, 517$, $583, 649, 671, 737, 781, 803, 869, 913, 979, 1067, 1111, 1133$.

The maximum distance after removing these residues is 14.

## 4. UPPER BOUNDS FOR THE MAXIMUM DISTANCE BETWEEN SUCCESSIVE PRIMES

Let's talk about the structure of RRS $p_r\#$. RRS $p_r\#$ includes RRS $p_{r-1}\#$ with the removal of residues that are multiples of $p_r$. Therefore, RRS $p_r\#$ consists of residues: $1, p_r, ..., p_{r-1}\# - p_r, p_{r-1}\# - 1, p_{r-1}\# + 1, p_{r-1}\# + p_r, ...$.

If one of the boundary values of RRS $p_{r-1}\#$, for example, $p_{r-1}\# - 1$ is a multiple of $p_r$, and it turns out to be a multiple $p_{r+1}$ of another boundary residue $p_{r-1}\# + 1$ in RRS $p_{r-1}\#$, then a distance equal to $p_{r-1}\# + p_r - (p_{r-1}\# - p_r) = 2p_r$ is formed between the residues of RRS $p_{r+1}\#$.

It was proved in [6] that the distance $2p_r$ exists in each RRS $p_{r+1}\#$, although it is not necessarily the largest.

Table A048670 QEIS shows that starting from RRS $2\#$ to RRS $19\#$, the maximum distance between RRS $p_r\#$ residues (Jakobsthal function value) is $j(p_r\#) = 2p_{r-1}$. Further, the values of the Jacobstal function $j(p_r\#) > 2p_{r-1}$ are encountered. The question arises - can there be a value of the Jacobstal function $j(p_r\#) < 2p_{r-1}$? It turns out that it doesn't. It was proved that $j(p_r\#) \geq 2p_{r-1}$ [6].

Let's now consider the RRS $p_r\#$ interval from $p_{r+1}$ to $p_{r+1}^2$, which contains only prime numbers. Let us denote the function of the maximum distance between successive numbers in RRS $p_r\#$ on the interval $p_{r+1}$ from to $p_{r+1}^2$ - $d(p_{r+1}^2)$. When passing from RRS $p_r\#$ to RRS $p_{r+1}\#$, the interval under consideration increases by $p_{r+2}^2 - p_{r+1}^2$, so the function $d(p_{r+1}^2)$ is monotonically increasing, but not strictly increasing, like the Jacobsthal function.

Conjecture 2

The following upper estimate is made for the maximum distance between successive prime numbers in RRS $p_r\#$ on the interval from $p_{r+1}$ to $p_{r+1}^2$:

$$d(p_{r+1}^2) \leq 2p_{r-1}. \qquad (4.1)$$

It was possible to calculate more values of the function (up to 4561#), due to the decrease in the interval when calculating the function $d(p_{r+1}^2)$ compared to the Jacobstal function.

The obtained values confirm the validity of estimate (4.1).

Table depending $d(p_{r+1}^2)$ on the value of the primorial

(only lines with a change in value are given)

| $p_r\#$ | $d(p_{r+1}^2)$ |
|---|---|
| 5# | 6 |
| 7# | 8 |
| 11# | 14 |
| 23# | 18 |
| 29# | 20 |
| 31# | 34 |
| 97# | 36 |
| 113# | 44 |
| 139# | 52 |
| 173# | 72 |
| 389# | 86 |
| 599# | 96 |
| 607# | 112 |
| 701# | 114 |
| 1153# | 118 |
| 1163# | 132 |
| 1409# | 148 |
| 2153# | 154 |
| 4129# | 180 |
| 4561# | 210 |

Taking into account that the inequality $j(p_r\#) \geq 2p_{r-1}$ was proved in [6], then if conjecture 2 is true and the interval from 1 to $p_{r+1}^2$ covers completely $p_r\#$ (more precisely, in view of the symmetry of RRS $p_r\#$, covers $p_r\#/2$), then the relation $\sup\{p_{r+1}-1, d(p_{r+1}^2)\} = 2p_{r-1}$ must hold. Let's prove it.

Assertion 3

If the condition $p_{r+1}^2 \geq p_r\#/2$ is true, then

$$\sup\{p_{r+1}-1, d(p_{r+1}^2)\} = 2p_{r-1}. \tag{4.2}$$

Proof

The relation $p_{r+1}^2 \geq p_r\#/2$ holds for $p_r = 7$, but the relation $p_{r+1}^2 \geq p_r\#/2$ does not hold for $p_r = 11$. Therefore, this ratio is also not satisfied for all $p_r > 11$ due to a faster increase $p_r\#/2$ than $p_{r+1}^2$.

Therefore, it is necessary to check the fulfillment of (4.2) for the values $p_r \leq 7$. We will use the table of dependence $d(p_{r+1}^2)$ on the value of the primorial when calculating.

The value is $\sup\{p_{r+1}-1, d(p_{r+1}^2)\} = \sup\{10, 8\} = 10 = 2 \times 5 = 2p_{r-1}$ for $p_r = 7$.

The value is $\sup\{p_{r+1}-1, d(p_{r+1}^2)\} = \sup\{6, 6\} = 6 = 2 \times 3 = 2p_{r-1}$ for $p_r = 5$.

The value is $\sup\{p_{r+1}-1, d(p_{r+1}^2)\} = \sup\{4, 4\} = 4 = 2 \times 2 = 2p_{r-1}$ for $p_r = 3$.

This corresponds to (4.2).

Assertion 3 shows that conjecture 2 does not contradict assertion $j(p_r\#) \geq 2p_{r-1}$.

Assertion 4

If $d(p_{r+1}^2) \leq 2p_r,$ \hfill (4.3)

then Legendre's conjecture holds.

Proof

Let's consider:

$$(p_r+1)^2 - p_r^2 = 2p_r + 1 > 2p_r, \tag{4.4}$$

therefore, if (4.3) holds, then there is at least one prime number between $p_r^2$ and $(p_r+1)^2$.

If we take $p_r+1$ instead of the value $p_r$ in (4.4), then we get:

$$(p_r+2)^2 - (p_r+1)^2 = 2p_r + 3 > 2p_r,$$

therefore, if (4.3) holds, then there is at least one prime number between $(p_r+1)^2$ and $(p_r+2)^2$ and so on.

If we take $p_r-1$ instead of the value $p_r$ in (4.4), then we get:

$$(p_r)^2 - (p_r-1)^2 = 2p_r + 1 > 2p_r, \qquad (4.5)$$

therefore, if (4.3) holds, then there is at least one prime number between $(p_r-1)^2$ and $(p_r)^2$.

Now if we take $p_r = 2$ in (4.4) or (4.5), then we obtain the fulfillment of Legendre's conjecture.

Conjecture 2 states - $d(p_{r+1}^2) \leq 2p_{r-1}$, therefore, having in mind that the function $d(p_{r+1}^2)$ is monotonically increasing, we get $d(p_{r+1}^2) \leq 2p_{r-1} \leq 2p_r$. Thus, when conjecture 2 is fulfilled, Legendre's conjecture is fulfilled.

## 5. RELATIONSHIP BETWEEN ASYMPTOTIC LOWER BOUNDS FOR THE OF THE JACOBSTAL FUNCTION AND THE MAXIMYM DISTANCE BETWEEN SUCCESSIVE PRIMES

Let $P^-(n) \leq y$ is a natural number without small prime factors $\leq y$.

Let $\Phi(x, y) = \{n \leq x : P^-(n) \leq y\}$ is the set of natural numbers not exceeding $x$, which have no prime factors $\leq y$.

All the asymptotic lower estimates $\sup_{p_{k+1} \leq x}(p_{k+1} - p_k)$ for available in the literature are the result of an asymptotic lower estimate for the Jacobstal function.

According to Rankin's theorem, if:

$$j(x\#) \gg x \ln x \ln \ln \ln x / (\ln \ln x)^2, \qquad (5.1)$$

then the relation is true:

$$\sup_{p_{k+1} \leq x}(p_{k+1} - p_k) \geq c \ln x \ln \ln x \ln \ln \ln x / (\ln \ln \ln x)^2. \qquad (5.2)$$

It was specified later that the constant in (5.2) is $c = e^\gamma + o(1)$.

The proof of Rankin's theorem using the properties of $\Phi(x, y)$ is given in [10].

An improved asymptotic lower estimate for the Jacobstal function compared to estimate (5.1) was obtained in [11]:

$$j(x\#) \gg x \ln x \ln \ln \ln x / \ln \ln x. \tag{5.3}$$

Based on a method similar to Rankin's theorem, the most accurate at this moment asymptotic lower estimate for the maximum distance between successive primes is proved there:

$$\sup_{p_{k+1} \leq x}(p_{k+1} - p_k) \geq c \ln x \ln \ln x \ln \ln \ln \ln x / \ln \ln \ln x. \tag{5.4}$$

If we compare formulas (5.1), (5.2), and (5.3), (5.4), then the asymptotic lower bounds for the distance between successive primes are obtained by a simple change of variable from the asymptotic lower bounds for the Jacobstal function.

Let us return to assertion 1 and put $p = x$ in formula (2.1):

$$\frac{\pi(x\#)}{\varphi(x\#)} = \frac{(e^\gamma + o(1)) \ln x}{x}. \tag{5.5}$$

The asymptotic of the ratio of the number of primes not exceeding $x\#$ and the number of residues in RRS $x\#$ is inversely proportional to the ratio of the asymptotic of the distance between successive primes and the distance between successive residues in RRS $x\#$.

Let's pay attention, that passing from the lower asymptotic of the Jacobstal function to the lower asymptotic of the maximum distance between successive primes, the value $(e^\gamma + o(1)) \ln x$ is substituted into a formula like (5.2) instead of the value $x$ in a formula like (5.1). It coincides with relation (5.5).

## 6. CONCLUSION AND SUGGESTIONS FOR FURTHER WORK

Work on estimating the maximum distance between successive primes is important and should be continued.

## 7. ACKNOWLEDGEMENTS

Thanks to everyone who has contributed to the discussion of this paper.